\documentclass[11pt]{article}
\usepackage{mathrsfs}
\usepackage{amssymb}
\usepackage{amsmath}
\usepackage[all]{xy}
\def\Im{\mathop{\rm Im}\nolimits}
\def\Ker{\mathop{\rm Ker}\nolimits}
\def\Coker{\mathop{\rm Coker}\nolimits}

\def\mod{\mathop{\rm mod}\nolimits}
\def\Mod{\mathop{\rm Mod}\nolimits}
\def\Hom{\mathop{\rm Hom}\nolimits}
\def\Ext{\mathop{\rm Ext}\nolimits}
\def\Tor{\mathop{\rm Tor}\nolimits}
\def\diag{\mathop{\rm diag}\nolimits}

\def\cx{\mathop{\rm cx}\nolimits}
\def\inf{\mathop{\rm inf}\nolimits}

\def\dim{\mathop{\rm dim}\nolimits}

\title{\Large \bf $n$-Strongly Gorenstein Projective, Injective and Flat Modules
}
\author{Guoqiang Zhao, Zhaoyong Huang\thanks{{\it E-mail address}:
zhaoguoqiang82@163.com, huangzy@nju.edu.cn}\\
{\footnotesize \it Department of Mathematics, Nanjing University,
Nanjing 210093, Jiangsu Province, P.R. China}}
\date{ }
\begin{document}
\baselineskip=18pt
\maketitle
\begin{abstract}
In this paper, we study the relation between $m$-strongly Gorenstein
projective (resp. injective) modules and $n$-strongly Gorenstein
projective (resp. injective) modules whenever $m \neq n$, and the
homological behavior of $n$-strongly Gorenstein projective (resp.
injective) modules. We introduce the notion of $n$-strongly
Gorenstein flat modules. Then we study the homological behavior of
$n$-strongly Gorenstein flat modules, and the relation between
these modules and $n$-strongly Gorenstein
projective (resp. injective) modules.
\end{abstract}

\vspace{0.5cm}

\centerline{\large \bf 1. Introduction}

\vspace{0.2cm}

As a nice generalization of the notion of finitely generated
projective modules, Auslander and Bridger introduced in [1] the
notion of finitely generated modules having Gorenstein dimension
zero over left and right Noetherian rings. For any module over a
general ring, Enochs and Jenda introduced in [8] the notion of
Gorenstein projective modules, which coincides with that of modules
having Gorenstein dimension zero for finitely generated modules over
left and right Noetherian rings. In [8] Enochs and Jenda also
introduced the dual notion of Gorenstein projective modules, which
is called Gorenstein injective modules. As a generalization of the
notion of flat modules, Enochs, Jenda and Torrecillas introduced in
[10] the notion of Gorenstein flat modules. These modules have been
studied extensively by many authors
(see [1, 2, 6, 7, 8, 9, 10, 15, 17], and so on). In particular, it
was proved that these modules share many nice
properties of the classical modules: projective,
injective and flat modules, respectively.

In 2007, Bennis and Mahdou introduced in [3] the notion of strongly
Gorenstein projective, injective, flat modules, which situate
between projective, injective, flat modules and Gorenstein
projective, injective, flat modules, respectively. Then they proved
that a module is Gorenstein projective (resp. injective) if and only
if it is a direct summand of a strongly Gorenstein projective (resp.
injective) module, and that every Gorenstein flat module is a direct
summand of a strongly Gorenstein flat module. Yang and Liu proved in
[18] that a module $M$ is strongly Gorenstein projective (resp.
injective, flat) if and only if so is $M\oplus H$ for any projective
(resp. injective, flat) module $H$. Gao and Zhang gave in [13] a
concrete construction of strongly Gorenstein projective modules, via
the existed construction of upper triangular matrix Artinian
algebras of degree two.

In a recent paper [4], for any $n\geq 1$, Bennis and Mahdou
introduced the notion of $n$-strongly Gorenstein projective and
injective modules, in which 1-strongly Gorenstein projective (resp.
injective) modules are just strongly Gorenstein projective (resp.
injective) modules. Then they proved that an $n$-strongly Gorenstein
projective module is projective if and only if it has finite flat
dimension. They also gave some equivalent characterizations of
$n$-strongly Gorenstein projective modules in terms of the vanishing
of some homological groups.

In this paper, based on the results mentioned above, we mainly study
the homological behavior of $n$-strongly Gorenstein projective,
injective and flat modules, and investigate the relation among them.
This paper is organized as follows.

In Section 2, we give the definitions of (strongly) Gorenstein
projective, injective and flat modules.

In Section 3, we study the relation between $m$-strongly Gorenstein
projective modules and $n$-strongly Gorenstein
projective modules whenever $m\neq n$, and the
closure of some special direct summand of an $n$-strongly Gorenstein
projective module. For any $n\geq 1$, we give an
example of an $n$-strongly Gorenstein projective
module, which is not $m$-strongly Gorenstein projective
whenever $n\nmid m$. For any $m, n \geq 1$, we prove that
the intersection of the subcategory of $m$-strongly Gorenstein
projective modules and that of $n$-strongly
Gorenstein projective modules is the subcategory
of $(m,n)$-strongly Gorenstein projective modules,
where $(m,n)$ is the greatest common divisor of $m$ and $n$. We give
a method how to construct a 1-strongly Gorenstein projective
module from $n$-strongly Gorenstein projective modules.
In addition, we prove that a module $M$ is
$n$-strongly Gorenstein projective if and only if
so is $M\oplus H$ for any projective module $H$,
which is a generalization of [18, Theorem 2.1]. We remark that
all the dual results hold for $n$-strongly Gorenstein
injective modules.

In Section 4, for any $n\geq 1$, we introduce the notion of
$n$-strongly Gorenstein flat modules, and then give an example of an
$n$-strongly Gorenstein flat module, which is not $m$-strongly
Gorenstein flat whenever $n\nmid m$. We prove that a module $M$ is
$n$-strongly Gorenstein flat if and only if so is $M\oplus H$ for
any flat module $H$. We also investigate the relation between
$n$-strongly Gorenstein flat modules and $n$-strongly Gorenstein
projective (resp.  injective) modules. We prove that a finitely
generated $n$-strongly Gorenstein projective module is finitely
presented $n$-strongly Gorenstein flat. In addition, we prove that
the character module of an $n$-strongly Gorenstein flat module is
$n$-strongly Gorenstein injective; and that the character module of
an $n$-strongly Gorenstein injective module is $n$-strongly
Gorenstein flat over an Artinian algebra. These results generalize
some results in [18].

\vspace{0.5cm}

\centerline{\large\bf 2. Preliminaries}

\vspace{0.2cm}

Throughout this paper, $R$ is an associative ring with identity and
$\Mod R$ is the category of left $R$-modules.

\vspace{0.2cm}

{\bf Definition 2.1} ([8]) A module $G \in \Mod R$ is called {\it
Gorenstein projective} ({\it G-projective} for short) if there
exists an exact sequence:
$$\cdots \to P_1 \to P_0 \to P^0 \to P^1\to \cdots ,$$
in $\Mod R$, such that: (1) All $P_i$ and $P^i$ are projective; (2)
$\Hom _{R}(-, P)$ leaves the sequence exact whenever $P\in \Mod R$
is projective; and (3) $G\cong \Im(P_0\to P^0)$. Dually, the notion of
{\it Gorenstein injective modules}
({\it G-injective modules} for short) is defined.

\vspace{0.2cm}

{\bf Definition 2.2} ([10]) A module $F \in \Mod R$ is called {\it
Gorenstein flat} ({\it G-flat} for short) if there exists an exact
sequence:
$$\cdots \to F_1 \to F_0 \to F^0 \to F^1\to \cdots ,$$
in $\Mod R$, such that: (1) All $F_i$ and $F^i$ are flat; (2)
$I\otimes _R-$ leaves the sequence exact whenever $I\in \Mod R^{op}$
is injective; and (3) $F\cong \Im(F_0\to F^0)$.

\vspace{0.2cm}

For a module $M\in \Mod R$, we denote $M^+=\Hom _{\mathbb{Z}}(M,
\mathbb{Q}/\mathbb{Z})$, which is called the {\it character module}
of $M$, where $\mathbb{Z}$ is the additive group of integers and
$\mathbb{Q}$ is the additive group of rational numbers (see [12]).
The following result is an analog of [15, Proposition 2.27], which is maybe
known.

\vspace{0.2cm}

{\bf Lemma 2.3} {\it A G-flat module with finite flat dimension is
flat.}

\vspace{0.2cm}

{\it Proof.} Let $M\in \Mod R$ be a G-flat module with finite flat
dimension. Then by [15, Theorem 3.6] and [12, Theorem 2.1], we have that $M^+\in \Mod R^{op}$ is
G-injective with finite injective dimension. So $M^+$ is injective
by the dual version of [15, Proposition 2.27], and hence $M$ is flat by [12, Theorem 2.1].
$\hfill{\square}$

\vspace{0.1cm}

{\bf Definition 2.4} ([3]) (1) A module $M\in \Mod R$ is called {\it
strongly Gorenstein projective} ({\it SG-projective} for short), if
there exists an exact sequence: $$0\rightarrow M\rightarrow P_0
\rightarrow M \rightarrow 0$$ in $\Mod R$ with $P_0$ projective,
such that $\Hom _{R}(-, P)$ leaves the sequence exact whenever $P\in
\Mod R$ is projective. Dually, the notion of
{\it strongly Gorenstein injective modules}
({\it SG-injective modules} for short) is defined.

(2) A module $M\in \Mod R$ is called {\it strongly Gorenstein flat}
({\it SG-flat} for short), if there exists an exact sequence:
$$0\rightarrow M\rightarrow F_0
\rightarrow M \rightarrow 0$$ in $\Mod R$ with $F_{0}$ flat, such
that $I\otimes_{R} -$ leaves the sequence exact whenever $I\in \Mod
R^{op}$ is injective.

\vspace{0.1cm}

It is trivial that $\{$projective modules$\}\subseteq
\{$SG-projective modules$\}\subseteq\{$G-projective modules$\}$,
$\{$injective modules$\}\subseteq \{$SG-injective
modules$\}\subseteq\{$G-injective modules$\}$ and $\{$flat
modules$\}\subseteq \{$SG-flat modules$\}\subseteq\{$G-flat
modules$\}$. By [3], all of the inclusions are strict in general.

\vspace{0.5cm}

\centerline{\large\bf 3. $n$-Strongly Gorenstein projective and
injective modules}

\vspace{0.2cm}

In this section we study the properties of $n$-strongly
Gorenstein projective modules. All the dual results hold for
the $n$-strongly Gorenstein injective
modules, and we omit this dual part.

\vspace{0.2cm}

{\bf Definition 3.1} ([4]) Let $n$ be a positive integer. A module
$M\in \Mod R$ is called {\it $n$-strongly Gorenstein projective}
({\it $n$-SG-projective} for short), if there exists an exact
sequence: $$0\rightarrow M\stackrel{f_{n}}{\longrightarrow} P_{n-1}
\stackrel{f_{n-1}}{\longrightarrow}
\cdots\stackrel{f_{1}}{\longrightarrow}P_{0}\stackrel{f_{0}}{\longrightarrow}M
\rightarrow 0$$ in $\Mod R$ with $P_{i}$ projective for any $0 \leq
i \leq n-1$, such that $\Hom _{R}(-, P)$ leaves the sequence exact
whenever $P\in \Mod R$ is projective. Dually, the notion of
{\it $n$-strongly Gorenstein injective modules}
({\it $n$-SG-injective modules} for short) is defined.

\vspace{0.2cm}

It is clear that the global dimension of $R$ is infinite if there
exists a non-projective $n$-SG-projective $R$-module for some $n \geq 1$.

In the following, we first study the relation between
$m$-SG-projective modules and $n$-SG-projective
modules whenever $m\neq n$.

Note that 1-SG-projective modules are just
SG-projective modules. In addition, for any $1
\leq i \leq n$, $\Im f_i$ in the above exact sequence
is also $n$-SG-projective. It is
trivial that a 1-SG-projective module is
$n$-SG-projective for any $n\geq 1$. However, for
any $n\geq 2$, an $n$-SG-projective module is not
necessarily $m$-SG-projective whenever $n\nmid m$,
as showed in the following example.

\vspace{0.2cm}

{\bf Example 3.2} Let $R$ be a finite-dimensional algebra over a
field given by the quiver:
$$
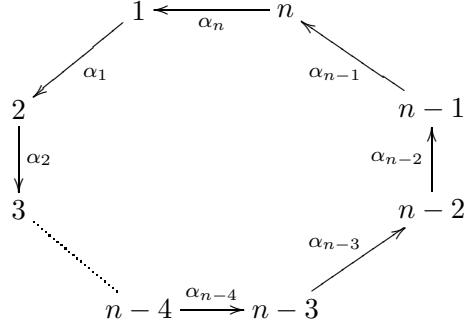
\xymatrix{
&1\ar[ld]^{\alpha_1}&n\ar[l]^{\alpha_{n}}&\\
2\ar[d]^{\alpha_2}&&&n-1\ar[ul]^{\alpha_{n-1}}\\
3&&&n-2\ar[u]^{\alpha_{n-2}}\\
&n-4\ar@{.}[lu]\ar[r]^{\alpha_{n-4}}&n-3\ar[ur]^{\alpha_{n-3}}& }
$$
modulo the ideal generated by
$\{\alpha_{i+1}\alpha_{i},\alpha_{1}\alpha_{n}$ $|$ $1\leq i\leq
n-1\}$. For any $1\leq i\leq n$, we use $S_i$, $P_i$ and $I^i$ to
denote the simple $R$-module, the indecomposable projective
$R$-module and the indecomposable injective $R$-module corresponding
to the vertex $i$, respectively. Then $R$ is a self-injective
algebra with infinite global dimension, and $P_{n}=I^1$,
$P_{i}=I^{i+1}$ for any $1\leq i\leq n-1$. In addition, for any
$1\leq i \leq n$, we have

(1) The following exact sequence
$$0\to S_i \to P_{i-1} \to \cdots \to P_1 \to P_n \to P_{n-1}
\to \cdots \to P_i \to S_i \to 0$$ is a minimal projective
resolution of $S_i$.

(2) For any $m\geq 1$, if $n\nmid m$, then $\Ext
^{m}_{R}(S_i,S_i)=0$; if $n\mid m$, then $\Ext
^{m}_{R}(S_i,S_i)\neq0$.

(3) $S_i$ is $n$-SG-projective.

(4) $S_i$ is not $m$-SG-projective whenever
$n\nmid m$.

\vspace{0.2cm}

For any $n\geq 1$, we use $n$-SG-Proj$(R)$
to denote the subcategory of $\Mod R$ consisting of
$n$-SG-projective modules. In the following,
assume that $m$ and $n$ are positive integers with $n\leq m$.

\vspace{0.2cm}

{\bf Lemma 3.3} {\it If $n\ |\ m$, then $n$-SG-Proj$(R)\subseteq
m$-SG-Proj$(R)$.}

\vspace{0.2cm}

We state a crucial result as follows.

\vspace{0.2cm}

{\bf Proposition 3.4} {\it (1) If $n\ |\ m$, then
$m$-SG-Proj$(R)\bigcap n$-SG-Proj$(R)=n$-SG-Proj$(R)$.

(2) If $n \nmid m$ and $m=kn+j$, where $k$ is a positive integer and
$0<j<n$, then $m$-SG-Proj$(R)\bigcap n$-SG-Proj$(R)\subseteq
j$-SG-Proj$(R)$. }

\vspace{0.2cm}

{\it Proof.} (1) It is trivial by Lemma 3.3.

(2) By Lemma 3.3, we have that $m$-SG-Proj$(R)\bigcap
n$-SG-Proj$(R)\subseteq m$-SG-Proj$(R)\bigcap kn$-SG-Proj$(R)$.
Assume that $M\in m$-SG-Proj$(R)\bigcap kn$-SG-Proj$(R)$. Then there
exists an exact sequence: $$0\rightarrow M \rightarrow
P_{m-1}\rightarrow\cdots\rightarrow P_1\rightarrow P_0\rightarrow
M\rightarrow 0\eqno{(1)}$$ in $\Mod R$ with $P_i$ projective for any
$0\leq i\leq m-1$. Put $L_{i}= \Ker(P_{i-1}\rightarrow P_{i-2})$ for
any $2 \leq i \leq m$. Because $M\in kn$-SG-Proj$(R)$, it is easy to
see that $M$ and $L_{kn}$ are projectively equivalent, that is,
there exist projective modules $P$ and $Q$ in $\Mod R$, such that
$M\oplus P$ $\cong$ $Q\oplus L_{kn}$.

First, consider the following pull-back diagram:
$$\xymatrix@C=20pt@R=15pt  {
    &&0\ar[d]&0\ar[d]\\
    &&Q\ar[d]\ar@{=}[r]&Q\ar[d]\\
    0\ar[r]&L_{kn+1}\ar[r]\ar@{=}[d]&X\ar[d]\ar[r]&M\oplus P\ar[r]\ar[d]&0\\
    0\ar[r]&L_{kn+1}\ar[r]&P_{kn}\ar[d]\ar[r]&L_{kn}\ar[r]\ar[d]&0\\
    &&0&0
} $$ Then $X$ is projective. Next, consider the following pull-back
diagram:
$$\xymatrix@=18pt  {
    &&0\ar[d]&0\ar[d]\\
    0\ar[r]&L_{kn+1}\ar[r]\ar@{=}[d]&Y\ar[d]\ar[r]&M\ar[r]\ar[d]&0\\
    0\ar[r]&L_{kn+1}\ar[r]&X\ar[d]\ar[r]&M\oplus P\ar[d]\ar[r]&0\\
    &&P\ar[d]\ar@{=}[r]&P\ar[d]\\
    &&0&0
    } $$
Thus $Y$ is also projective. Combining the exact sequence (1) and
the first row in the above diagram, we get the following exact
sequence: $$0\rightarrow M \rightarrow
P_{m-1}\rightarrow\cdots\rightarrow P_{kn+1}\rightarrow Y\rightarrow
M\rightarrow 0,$$ which is still exact after applying the functor
$\Hom _R(-,P)$ for any projective $R$-module $P$. So $M$ is
$j$-SG-projective, and hence $m$-SG-Proj$(R)\bigcap
n$-SG-Proj$(R)\subseteq j$-SG-Proj$(R)$.  $\hfill{\square}$

\vspace{0.2cm}

We use $(m, n)$ to denote the greatest common divisor of $m$ and
$n$.

\vspace{0.2cm}

{\bf Theorem 3.5} {\it $m$-SG-Proj$(R)\bigcap
n$-SG-Proj$(R)=(m,n)$-SG-Proj$(R)$.}

\vspace{0.2cm}

{\it Proof.} If $n\ |\ m$, then the assertion follows from
Proposition 3.4(1).

Now assume that $n \nmid m$ and $m=k_{0}n+j_0$, where $k_0$ is a
positive integer and $0<j_{0}<n$. By Proposition 3.4(2), we have
$m$-SG-Proj$(R)\bigcap n$-SG-Proj$(R)\subseteq j_0$-SG-Proj$(R)$. If
$j_0 \nmid n$ and $n=k_{1}j_{0}+j_1$ with $0<j_{1}<j_0$, then by
Proposition 3.4(2) again, we have that $m$-SG-Proj$(R)\bigcap
n$-SG-Proj$(R)\subseteq n$-SG-Proj$(R)\bigcap
j_0$-SG-Proj$(R)\subseteq j_1$-SG-Proj$(R)$. Continuing the above
procedure, after finite steps, there exists a positive integer $t$
such that $j_t=k_{t+2}j_{t+1}$ and $j_{t+1}=(m, n)$. Thus
$m$-SG-Proj$(R)\bigcap n$-SG-Proj$(R)\subseteq
j_t$-SG-Proj$(R)\bigcap
j_{t+1}$-SG-Proj$(R)=j_{t+1}$-SG-Proj$(R)=(m, n)$-SG-Proj$(R)$. On
the other hand, we always have $(m, n)$-SG-Proj$(R) \subseteq
m$-SG-Proj$(R)\bigcap n$-SG-Proj$(R)$, so they are identical.
$\hfill{\square}$

\vspace{0.2cm}

As an immediate consequence of Theorem 3.5, we have the following

\vspace{0.2cm}

{\bf Corollary 3.6} {\it $n$-SG-Proj$(R)\bigcap (n+1)$-SG-Proj$(R) =
$1-SG-Proj$(R)$. In particular,
$\bigcap_{n\geq2}n$-SG-Proj$(R)=$1-SG-Proj$(R)$.}

\vspace{0.2cm}

The following result shows that the difference between the
projectivity and $n$-SG-projectivity of modules is
the self-orthogonality of modules.

\vspace{0.2cm}

{\bf Proposition 3.7} {\it Let $M\in \Mod R$ be $n$-SG-projective
and $n\geq 1$. The following statements are
equivalent.

(1) $M$ is projective.

(2) $\Ext ^{i}_{R}(M, M)=0$ for any $i \geq 1$.

(3) $\Ext ^{i}_{R}(M, M)=0$ for any $1\leq i\leq n$.}

\vspace{0.2cm}

{\it Proof.} $(1)\Rightarrow (2)\Rightarrow (3)$ are trivial. By the
dimension shifting, it is easy to get $(3)\Rightarrow (1)$.
$\hfill{\square}$

\vspace{0.2cm}

In the rest of this section, we will study the homological behavior
of $n$-SG-projective modules.

\vspace{0.2cm}

{\bf Proposition 3.8} {\it For any $n\geq 1$,
$n$-SG-Proj$(R)$ is closed under direct sums.}

\vspace{0.2cm}

{\it Proof.} Let $\{M_j\}_{j\in J}$ be a family of $n$-SG-projective
modules in $\Mod R$. Then for any $j \in J$, there exists an exact
sequence: $$0\rightarrow M_j\rightarrow P^{(j)}_{n-1} \rightarrow
\cdots\rightarrow P^{(j)}_{0}\rightarrow M_j \rightarrow 0$$ in
$\Mod R$ with $P^{(j)}_{i}$ projective for any $0 \leq i \leq n-1$,
such that $\Hom _{R}(-, P)$ leaves the sequence exact whenever $P\in
\Mod R$ is projective. So we get an exact sequence: $$0\rightarrow
\oplus _{j\in J}M_j\rightarrow \oplus _{j\in J}P^{(j)}_{n-1}
\rightarrow \cdots\rightarrow \oplus _{j\in J}P^{(j)}_{0}\rightarrow
\oplus _{j\in J}M_j \rightarrow 0$$ in $\Mod R$. Because $\oplus
_{j\in J}P^{(j)}_{n-1}, \cdots, \oplus _{j\in J}P^{(j)}_{0}$ are
projective and the obtained exact sequence is still exact after
applying the functor $\Hom _{R}(-, P)$ whenever $P\in \Mod R$ is
projective, $\oplus _{j\in J}M_j$ is $n$-SG-projective and the
assertion follows. $\hfill{\square}$

\vspace{0.2cm}

The following result gives some characterizations of
$n$-SG-projective modules, which also gives a method how to
construct a 1-SG-projective module from $n$-SG-projective modules.

\vspace{0.2cm}

{\bf Theorem 3.9} {\it For any $M\in \Mod R$ and $n\geq 1$, the
following statements are equivalent.

(1) $M$ is $n$-SG-projective.

(2) There exists an exact sequence:
$$0\rightarrow M\stackrel{f_{n}}{\longrightarrow} P_{n-1}
\stackrel{f_{n-1}}{\longrightarrow}
\cdots\stackrel{f_{1}}{\longrightarrow}P_{0}\stackrel{f_{0}}{\longrightarrow}M
\rightarrow 0$$ in $\Mod R$ with $P_{i}$ projective for any $0\leq
i\leq n-1$, such that $\oplus_{i=1}^{n}\Im f_{i}$ is
1-SG-projective.

(3) There exists an exact sequence:
$$0\rightarrow M\stackrel{f_{n}}{\longrightarrow} P_{n-1}
\stackrel{f_{n-1}}{\longrightarrow}
\cdots\stackrel{f_{1}}{\longrightarrow}P_{0}\stackrel{f_{0}}{\longrightarrow}M
\rightarrow 0$$ in $\Mod R$ with $P_{i}$ projective for any $0\leq
i\leq n-1$, such that $\oplus_{i=1}^{n}\Im f_{i}$ is G-projective.

(4) There exists an exact sequence:
$$0\rightarrow M\stackrel{f_{n}}{\longrightarrow} P_{n-1}
\stackrel{f_{n-1}}{\longrightarrow}
\cdots\stackrel{f_{1}}{\longrightarrow}P_{0}\stackrel{f_{0}}{\longrightarrow}M
\rightarrow 0$$ in $\Mod R$, where $P_{i}$ has finite projective
dimension for any $0\leq i\leq n-1$, such that $\oplus_{i=1}^{n}\Im
f_{i}$ is 1-SG-projective.

(5) There exists an exact sequence:
$$0\rightarrow M\stackrel{f_{n}}{\longrightarrow} P_{n-1}
\stackrel{f_{n-1}}{\longrightarrow}
\cdots\stackrel{f_{1}}{\longrightarrow}P_{0}\stackrel{f_{0}}{\longrightarrow}M
\rightarrow 0$$ in $\Mod R$, where $P_{i}$ has finite projective
dimension for any $0\leq i\leq n-1$, such that $\oplus_{i=1}^{n}\Im
f_{i}$ is G-projective.}

\vspace{0.2cm}

{\it Proof.} $(1)\Rightarrow (2)$ Let $M\in \Mod R$ be
$n$-SG-projective. Then we have an exact sequence:
$$0\rightarrow M\stackrel{f_{n}}{\longrightarrow} P_{n-1}
\stackrel{f_{n-1}}{\longrightarrow}
\cdots\stackrel{f_{1}}{\longrightarrow}P_{0}\stackrel{f_{0}}{\longrightarrow}M
\rightarrow 0$$ in $\Mod R$ with $P_{i}$ projective for any $0 \leq
i \leq n-1$, such that $\Hom _{R}(-, P)$ leaves the sequence exact
whenever $P\in \Mod R$ is projective. Thus, for each $1\leq i\leq
n$, we have an exact sequence:
$$0\rightarrow \Im f_i\stackrel{\alpha _{i}}{\longrightarrow}
P_{i-1}\stackrel{f_{i-1}}{\longrightarrow}\cdots
\stackrel{f_{1}}{\longrightarrow}
P_{0}\stackrel{f_{n}}{\longrightarrow}
P_{n-1}\stackrel{f_{n-1}}{\longrightarrow}\cdots\stackrel{f_{i+1}}{\longrightarrow}
P_i\stackrel{f_{i}}{\longrightarrow} \Im f_i\rightarrow 0$$ in $\Mod
R$. By adding these exact sequences, we get the following exact
sequence:
$$0\rightarrow \oplus _{i=1}^n\Im f_i
\stackrel{\alpha}{\longrightarrow} \oplus
_{i=0}^{n-1}P_i\stackrel{f}{\longrightarrow} P_{n-1}\oplus
P_{0}\oplus\cdots\oplus P_{n-2}\rightarrow \cdots,$$ where $\alpha
=\diag\{\alpha _1, \alpha _2, \cdots, \alpha _n\}$ and $f
=\diag\{f_nf_0, f_1, \cdots, f_{n-1}\}$. It is easy to see that $\Im
f\cong \oplus _{i=1}^n\Im f_i$ and $\Ext _R^1(\oplus _{i=1}^n\Im
f_i, P)=0$ for any projective module $P\in \Mod R$, which implies
$\oplus_{i=1}^{n} \Im f_{i}$ is 1-SG-projective.

$(2)\Rightarrow (3)\Rightarrow (5)$ and $(2)\Rightarrow
(4)\Rightarrow (5)$ are trivial.

$(5)\Rightarrow (1)$ Assume that
$$0\rightarrow M\stackrel{f_{n}}{\longrightarrow} P_{n-1}
\stackrel{f_{n-1}}{\longrightarrow}
\cdots\stackrel{f_{1}}{\longrightarrow}P_{0}\stackrel{f_{0}}{\longrightarrow}M
\rightarrow 0$$ is an exact sequence in $\Mod R$, where $P_{i}$ has
finite projective dimension for any $0\leq i\leq n-1$, such that
$\oplus_{i=1}^{n}\Im f_{i}$ is G-projective. Then, for any $0\leq
i\leq n-1$, we have the exact sequence: $$0\rightarrow \Im
f_{i+1}\rightarrow P_{i}\rightarrow \Im f_{i}\rightarrow 0.$$
Because $\oplus_{i=1}^{n} \Im f_{i}$ is G-projective, so is each
$P_i$ by [15, Theorem 2.5]. Thus each $P_i$ is projective by [15, Proposition 2.27].
In particular, $M$ is also G-projective by [15, Theorem 2.5], so
$\Ext ^{i}_{R}(M, P)$ = 0 for any projective module $P\in \Mod R$
and $i\geq 1$. It follows that $M$ is $n$-SG-projective.
$\hfill{\square}$

\vspace{0.2cm}

From [18] we know that 1-SG-Proj$(R)$ is not
closed under direct summands. The following example illustrates that
for any $n \geq 1$, $n$-SG-Proj$(R)$ is not
closed under direct summands.

\vspace{0.2cm}

{\bf Example 3.10} Under the assumption of Example 3.2, we have
that $\oplus_{i=1}^{n}S_{i}$ is 1-SG-projective by Theorem 3.9, and hence
$(n-1)$-SG-projective. However, for
any $1\leq i \leq n$, $S_i$ is not $(n-1)$-SG-projective.

\vspace{0.2cm}

The following result is a generalization of [18, Theorem 2.1], which
shows that some special direct summand of an $n$-SG-projective
module is again $n$-SG-projective. For a module $M\in \Mod R$, we
use $\underline{M}$ to denote the maximal submodule of $M$ without
projective summands.

\vspace{0.2cm}

{\bf Theorem 3.11} {\it For any $n\geq 1$, a module $M\in \Mod R$ is
$n$-SG-projective if and only if so is $\underline{M}$.}

\vspace{0.2cm}

{\it Proof.} Let $M=\underline{M}\oplus P$ with $P$ a projective
module in $\Mod R$. If $\underline{M}$ is $n$-SG-projective, then
$M$ is also $n$-SG-projective by Proposition 3.8.

Conversely, assume that $M\in \Mod R$ is $n$-SG-projective. Then
there exists an exact sequence: $$0\rightarrow
(M=)\underline{M}\oplus P\stackrel{f_{n}}{\longrightarrow} P_{n-1}
\stackrel{f_{n-1}}{\longrightarrow}
\cdots\stackrel{f_{1}}{\longrightarrow}P_{0}
\stackrel{f_{0}}{\longrightarrow}\underline{M}\oplus
P(=M)\rightarrow 0$$ in $\Mod R$ with $P_{i}$ projective for any
$0\leq i \leq n-1$, such that $\Hom _{R}(-,Q)$ leaves the sequence
exact whenever $Q\in \Mod R$ is projective.

Put $\Im f_{i}=K_{i}$ for any $0\leq i \leq n$. First, consider the
following push-out diagram:
$$\xymatrix@=18pt  {
&&0\ar[d]&0\ar[d]\\
0\ar[r]&P\ar[r]\ar@{=}[d]&M\ar[d]^{f_n}\ar[r]&\underline{M}\ar[r]\ar[d]&0\\
0\ar[r]&P\ar[r]&P_{n-1}\ar[d]\ar[r]&Q_{n-1}\ar[d]\ar[r]&0\\
&&K_{n-1}\ar[d]\ar@{=}[r]&K_{n-1}\ar[d]\\
&&0&0 } $$ Because $M$ is G-projective, both $\underline{M}$ and $Q_{n-1}$ are also
G-projective by [15, Theorem 2.5]. So Ext$^{1}_{R} (Q_{n-1}, P)=0$ and
the middle row $0\rightarrow P\rightarrow P_{n-1}\rightarrow
Q_{n-1}\rightarrow 0$ in the above diagram splits, which implies
that $Q_{n-1}$ is projective. Because $K_{n-1}$ is also
$n$-SG-projective, the third column
$$0\rightarrow \underline{M}\rightarrow Q_{n-1}\rightarrow
K_{n-1}\rightarrow 0$$ in the above diagram is still exact after
applying the functor $\Hom _{R}(-,Q)$ whenever $Q\in \Mod R$ is
projective.

Next, consider the following pull-back diagram:
$$\xymatrix@C=20pt@R=15pt {
&0\ar[d]&0\ar[d]\\
&K_1\ar[d]\ar@{=}[r]&K_1\ar[d]\\
0\ar[r]&Q_0\ar[r]\ar[d]&P_0\ar[d]^{f_0}\ar[r]&P\ar[r]\ar@{=}[d]&0\\
0\ar[r]&\underline{M}\ar[d]\ar[r]&M\ar[d]\ar[r]&P\ar[r]&0\\
&0&0}\ \ \ $$ Then $0\rightarrow K_1\rightarrow Q_0\rightarrow
\underline{M}\rightarrow 0$ is exact and $Q_0$ is projective. Thus
we obtain the following exact sequence:
$$0\rightarrow \underline{M}\rightarrow Q_{n-1}\rightarrow
P_{n-2}\rightarrow\cdots\rightarrow P_{1}\rightarrow Q_0\rightarrow
\underline{M}\rightarrow 0.$$ Note that both $K_1$ and
$\underline{M}$ are G-projective. Thus the above exact sequence is
still exact after applying the functor $\Hom _{R}(-,Q)$ whenever
$Q\in \Mod R$ is projective, which implies $\underline{M}$ is
$n$-SG-projective. $\hfill{\square}$

\vspace{0.2cm}

By Theorem 3.11, we immediately have the following

\vspace{0.2cm}

{\bf Corollary 3.12} {\it Assume that $M, N\in \Mod R$ are
projectively equivalent. Then, for
any $n\geq 1$, $M$ is $n$-SG-projective if and
only if so is $N$.}

\vspace{0.2cm}

We denote $\mod R$ the category of finitely generated left
$R$-modules, and $n$-SG-proj$(R)=\{M\in \mod R \ |$ there exists an
exact sequence $0\rightarrow M\rightarrow P_{n-1}\rightarrow
P_{n-2}\rightarrow\cdots\rightarrow P_0\rightarrow M\rightarrow 0$
in $\mod R$ with $P_i$ projective for any $0\leq i\leq n-1$, such
that $\Hom _{R}(-,P)$ leaves the sequence exact whenever
$P\in \mod R$ is projective$\}$.

The following fact is useful, which is a generalization of [13,
Proposition 1.1].

\vspace{0.2cm}

{\bf Lemma 3.13} {\it For any $n\geq 1$, $n$-SG-Proj$(R)\bigcap \mod
R = n$-SG-proj$(R)$.}

\vspace{0.2cm}

{\it Proof.} Let $M\in n$-SG-proj$(R)$. By using an argument similar
to that of [13, Proposition 1.1], we have that $M\in$
$n$-SG-Proj$(R)\bigcap\mod R$.

Conversely, let $M\in$ $n$-SG-Proj$(R)\bigcap\mod R$. Then there
exists an exact sequence:
$$0\rightarrow M\stackrel{f_{n}}{\longrightarrow} P_{n-1}
\stackrel{f_{n-1}}{\longrightarrow}
\cdots\stackrel{f_{1}}{\longrightarrow}P_{0}\stackrel{f_{0}}{\longrightarrow}M
\rightarrow 0\eqno{(2)}$$ in $\Mod R$ with $P_{i}$ projective for
any $0\leq i \leq n-1$, such that $\Hom _{R}(-,P)$ leaves the
sequence exact whenever $P\in \Mod R$ is projective. Put $\Im f_i=
K_i$ for any $0\leq i \leq n$. There exists a projective module
$P^{\prime}_{n-1}\in \Mod R$ such that $P_{n-1}\oplus
P^{\prime}_{n-1}=Q$ is free, so we have an exact sequence:
$$0\rightarrow
M\stackrel{f^{\prime}_{n}}{\longrightarrow}P_{n-1}\oplus
P^{\prime}_{n-1}\stackrel{f_{n-1}^\prime}{\longrightarrow}
P_{n-2}\oplus
P^{\prime}_{n-1}\stackrel{f_{n-2}^\prime}{\longrightarrow} P_{n-3}
\stackrel{f_{n-3}}{\longrightarrow}\cdots
\stackrel{f_{1}}{\longrightarrow}P_{0}\stackrel{f_{0}}{\longrightarrow}M
\rightarrow 0.$$ Then $\Im f_{n-1}^\prime\cong K_{n-1}\oplus
P^{\prime}_{n-1}$ and $\Im f_{n-2}^\prime\cong K_{n-2}$. Since $M$
is finitely generated, one can write $Q = Q_{n-1}\oplus
Q^{\prime}_{n-1}$ with $Q_{n-1}\in \mod R$ and $\Im
f_{n}^\prime\subseteq Q_{n-1}$. So we get an exact sequence:
$$0\longrightarrow M\stackrel{f^{\prime}_{n}}{\longrightarrow}
Q_{n-1}\longrightarrow K^{\prime}_{n-1}\longrightarrow 0\eqno{(3)}$$
with $K^{\prime}_{n-1}\oplus Q^{\prime}_{n-1}\cong\Im
f_{n-1}^\prime$, and hence $K^{\prime}_{n-1}\in
n$-SG-Proj$(R)\bigcap\mod R$ by Corollary 3.12.

Consider the following push-out diagram:
$$\xymatrix@C=20pt@R=15pt {
&0\ar[d]&0\ar[d]\\
&Q^{\prime}_{n-1}\ar[d]\ar@{=}[r]&Q^{\prime}_{n-1}\ar[d]\\
0\ar[r]&\Im f_{n-1}^\prime\ar[r]\ar[d]&P_{n-2}\oplus
P^{\prime}_{n-1}\ar[d]\ar[r]&K_{n-2}\ar[r]\ar@{=}[d]&0\\
0\ar[r]&K^{\prime}_{n-1}\ar[d]\ar[r]&X\ar[d]\ar[r]&K_{n-2}\ar[r]&0\\
&0&0}\ \ \ $$ Then $X$ is G-projective by [15, Theorem 2.5], and so the
middle column $0\rightarrow Q^{\prime}_{n-1}\rightarrow
P_{n-2}\oplus P^{\prime}_{n-1}\rightarrow X\rightarrow 0$ in the
above diagram splits, which implies that $X$ is projective.
Combining the exact sequences (2), (3) with the third row in the
above diagram, we get an exact sequence:
$$0\rightarrow
M\stackrel{f^{\prime}_{n}}{\longrightarrow} Q_{n-1}\rightarrow
X\rightarrow P_{n-3}\rightarrow\cdots\rightarrow P_{0}\rightarrow
M\rightarrow 0$$ with $Q_{n-1}\in \mod R$. Repeating the above
procedure with $K^{\prime}_{n-1}(\cong\Coker f^{\prime}_{n})$
replacing $M$, we finally obtain the following exact sequence:
$$0\rightarrow M\rightarrow Q_{n-1}\rightarrow
Q_{n-2}\rightarrow\cdots\rightarrow Q_0\rightarrow M\rightarrow 0$$
in $\mod R$ with $Q_i$ projective for any $0\leq i\leq n-1$, which
implies $M\in n$-SG-proj$(R)$. $\hfill{\square}$

\vspace{0.2cm}

The following result gives some equivalent characterizations of
finitely generated $n$-SG-projective modules.

\vspace{0.2cm}

{\bf Theorem 3.14} {\it For any $M\in \mod R$ and $n\geq 1$, the
following statements are equivalent.

(1) M is $n$-SG-projective.

(2) There exists an exact sequence: $$0\rightarrow M\rightarrow
P_{n-1}\rightarrow P_{n-2}\rightarrow\cdots\rightarrow
P_0\rightarrow M\rightarrow 0$$ in $\mod R$ with each $P_i$
projective for any $0\leq i\leq n-1$, and $\Ext ^{i}_{R}(M, R) = 0$
for any $i\geq 1$.

(3) There exists an exact sequence: $$0\rightarrow M\rightarrow
P_{n-1}\rightarrow P_{n-2}\rightarrow\cdots\rightarrow
P_0\rightarrow M\rightarrow 0$$ in $\mod R$ with $P_i$ projective
for any $0\leq i\leq n-1$, and $\Ext ^{i}_{R}(M, F) = 0$ for any
flat module $F\in \Mod R$ and $i\geq 1$.

(4) There exists an exact sequence: $$0\rightarrow M\rightarrow
P_{n-1}\rightarrow P_{n-2}\rightarrow\cdots\rightarrow
P_0\rightarrow M\rightarrow 0$$ in $\mod R$ with $P_i$ projective
for any $0\leq i\leq n-1$, and $\Ext ^{i}_{R}(M, F) = 0$ for any
$F\in \Mod R$ with finite flat dimension and $i\geq 1$.}

\vspace{0.2cm}

{\it Proof.} $(1)\Rightarrow (2)$ follows from Lemma 3.13. The
proofs of other implications are similar to that of [3, Proposition
2.12], so we omit them. $\hfill{\square}$

\vspace{0.2cm}

We have obtained some properties of the intersection between
$m$-SG-projective modules and $n$-SG-projective modules (see 3.3--3.6).
We end this section with some properties of the union of these modules.

It has been known that $\bigcup _{n\geq 1}n$-SG-Proj$(R)
\subseteq\{$G-projective $R$-modules$\}$. We will
show that this inclusion is strict in general, and
also investigate when the equality holds true.

In the rest of this section, $R$ is a finite-dimensional $k$-algebra
over an algebraically closed field $k$. Let $M\in \mod R$ and
$$\cdots\rightarrow P_{n}\rightarrow\cdots\rightarrow
P_{1}\rightarrow P_{0}\rightarrow M\rightarrow 0$$ be a minimal
projective resolution of $M$ in $\mod R$. Recall from [11] that {\it
complexity} of $M$ is defined as $\cx(M)= \inf \{ b\geq 0\ |$ there
exists a $c>0$ such that $\dim _{k}P_{n}\leq cn^{b-1}$ for all $n\}$
if it exists, otherwise $\cx(M)=\infty$. It is easy to see that
$\cx(M)=0$ implies $M$ is of finite projective dimension, and
$\cx(M)\leq1$ if and only if the dimensions of $P_{n}$ are bounded.

\vspace{0.2cm}

{\bf Proposition 3.15} {\it Let $R$ be a self-injective algebra.

(1) If $R$ is of infinite representation type with vanishing radical
cube, then $\bigcup _{n\geq 1}n$-SG-Proj$(R)
\subsetneqq\{$G-projective $R$-modules$\}$.

(2) If $R$ is of finite representation type, then $\bigcup _{n\geq
1}n$-SG-proj$(R)=\{$finitely generated G-projective $R$-modules$\}$.}

\vspace{0.2cm}

{\it Proof.} Let $R$ be a self-injective algebra. Then
$\mod R =\{$finitely generated G-projective
$R$-modules$\}$.

(1) Assume that $R$ is of infinite representation type with
vanishing radical cube. Then by [14, Theorem 6.1], there exists a
module $M\in \mod R$ such that $\cx(M)\geq 2$. It is easy to see
that $M$ is not $n$-SG-projective for any $n \geq 1$. Thus $\bigcup
_{n\geq 1}n$-SG-proj$(R) \subsetneqq\{$finitely generated
G-projective $R$-modules$\}$, and therefore $\bigcup _{n\geq
1}n$-SG-Proj$(R) \subsetneqq\{$G-projective $R$-modules$\}$ by Lemma
3.13.

(2) Assume that $R$ is of finite representation type. We claim that
any indecomposable module $M\in \mod R$ is $n$-SG-projective for
some $n \geq 1$. Otherwise, if $M\in \mod R$ is not
$n$-SG-projective for any $n \geq 1$. Then there exists a minimal
projective resolution:
$$\cdots \to P_i \to \cdots \to P_1 \to P_0 \to M \to 0$$
of $M$ in $\mod R$, which is of infinite length. Because $R$ is
self-injective, all $P_i$ are also injective. Then by [16, Lemma 2.6],
all syzysy modules in the above exact sequence are indecomposable.
It is not difficult to see that any two of these syzysy modules are
not isomorphic, which implies that $R$ is of infinite representation
type. This is a contradiction. The claim is proved. So it follows
from Proposition 3.8 that any module $M\in \mod R$ is
$n$-SG-projective for some $n \geq 1$. Thus we get that $\bigcup
_{n\geq 1}n$-SG-proj$(R)=\{$finitely generated G-projective
$R$-modules$\}$. $\hfill{\square}$

\vspace{0.5cm}

\centerline{\large\bf 4. $n$-Strongly Gorenstein flat modules}

\vspace{0.2cm}

In this section, we introduce the notion of $n$-strongly Gorenstein
flat modules. Then we study the homological behavior of $n$-strongly
Gorenstein flat modules, and the relation between these modules and
$n$-strongly Gorenstein projective (resp. injective) modules.

\vspace{0.2cm}

{\bf Definition 4.1} Let $n$ be a positive integer. A module $M\in
\Mod R$ is called {\it $n$-strongly Gorenstein flat} ({\it
$n$-SG-flat} for short), if there exists an exact sequence:
$$0\rightarrow M\stackrel{h_{n}}{\longrightarrow} F_{n-1}\stackrel{h_{n-1}}{\longrightarrow}
F_{n-2}\stackrel{h_{n-2}}{\longrightarrow}\cdots\stackrel{h_{1}}{\longrightarrow}
F_0\stackrel{h_{0}}{\longrightarrow} M\rightarrow 0$$ in $\Mod R$
with $F_{i}$ flat for any $0\leq i \leq n-1$, such that
$I\otimes_{R} -$ leaves the sequence exact whenever $I\in \Mod
R^{op}$ is injective.

\vspace{0.2cm}

Note that 1-SG-flat modules are just SG-flat modules. For any $1\leq
i\leq n$, $\Im h_i$ in the above exact sequence is also $n$-SG-flat.

Let $n$ be a positive integer. It is trivial that a 1-SG-flat
(especially, flat) module is $n$-SG-flat, and an $n$-SG-flat module
is G-flat. It is clear that the weak global dimension of $R$ is
infinite if there exists a non-flat $n$-SG-flat $R$-module for some
$n \geq 1$. On the other hand, for a quasi-Frobenius ring $R$, it is
easy to see that a module in $\Mod R$ is $n$-SG-flat if and only if
it is $n$-SG-projective, if and only if it is $n$-SG-injective. So
we have the following example which illustrates that there exists an
$n$-SG-flat module, but it is not $m$-SG-flat whenever $n\nmid m$.

\vspace{0.2cm}

{\bf Example 4.2} Under the assumption of Example 3.2, because $R$
is quasi-Frobenius, for any $1 \leq i \leq n$, we have the following
facts: (1) $S_i$ is $n$-SG-flat; and (2) $S_i$ is not $m$-SG-flat
whenever $n\nmid m$.

\vspace{0.2cm}

{\bf Proposition 4.3} {\it For any $n \geq 1$, the subcategory
$n$-SG-Flat$(R)$ of $\Mod R$ consisting of $n$-SG-flat modules is
closed under direct sums.}

\vspace{0.2cm}

{\it Proof.} The proof is similar to that of Proposition 3.8, so we
omit it.  $\hfill{\square}$

\vspace{0.2cm}

The following result is an analog of Theorem 3.9, which gives some
characterizations of $n$-SG-flat modules, and also gives a method
how to construct a 1-SG-flat module from $n$-SG-flat modules.

\vspace{0.2cm}

{\bf Theorem 4.4} {\it For any $M\in \Mod R$ and $n\geq 1$, consider
the following conditions.

(1) $M$ is $n$-SG-flat.

(2) There exists an exact sequence:
$$0\rightarrow M\stackrel{h_{n}}{\longrightarrow} F_{n-1}
\stackrel{h_{n-1}}{\longrightarrow}
\cdots\stackrel{h_{1}}{\longrightarrow}F_{0}\stackrel{h_{0}}{\longrightarrow}M
\rightarrow 0$$ in $\Mod R$ with $F_{i}$ flat for any $0\leq i\leq
n-1$, such that $\oplus_{i=1}^{n}\Im h_{i}$ is 1-SG-flat.

(3) There exists an exact sequence:
$$0\rightarrow M\stackrel{h_{n}}{\longrightarrow} F_{n-1}
\stackrel{h_{n-1}}{\longrightarrow}
\cdots\stackrel{h_{1}}{\longrightarrow}F_{0}\stackrel{h_{0}}{\longrightarrow}M
\rightarrow 0$$ in $\Mod R$ with $F_{i}$ flat for any $0\leq i\leq
n-1$, such that $\oplus_{i=1}^{n}\Im h_{i}$ is G-flat.

(4) There exists an exact sequence:
$$0\rightarrow M\stackrel{h_{n}}{\longrightarrow} F_{n-1}
\stackrel{h_{n-1}}{\longrightarrow}
\cdots\stackrel{h_{1}}{\longrightarrow}F_{0}\stackrel{h_{0}}{\longrightarrow}M
\rightarrow 0$$ in $\Mod R$, where $F_{i}$ has finite flat dimension
for any $0\leq i\leq n-1$, such that $\oplus_{i=1}^{n}\Im h_{i}$ is
1-SG-flat.

(5) There exists an exact sequence:
$$0\rightarrow M\stackrel{h_{n}}{\longrightarrow} F_{n-1}
\stackrel{h_{n-1}}{\longrightarrow}
\cdots\stackrel{h_{1}}{\longrightarrow}F_{0}\stackrel{h_{0}}{\longrightarrow}M
\rightarrow 0$$ in $\Mod R$, where $F_{i}$ has finite flat dimension
for any $0\leq i\leq n-1$, such that $\oplus_{i=1}^{n}\Im h_{i}$ is
G-flat.

In general, we have $(1)\Leftrightarrow (2) \Leftrightarrow
(3)\Rightarrow (4)\Rightarrow (5)$. If $R$ is a right coherent ring,
then all of these conditions are equivalent.}

\vspace{0.2cm}

{\it Proof.} $(1)\Rightarrow (2)$ By using an argument similar to
that in the proof of $(1)\Rightarrow (2)$ in Theorem 3.9, we get
the assertion.

$(2)\Rightarrow (3)\Rightarrow (5)$ and $(2)\Rightarrow
(4)\Rightarrow (5)$ are trivial, and it is easy to get
$(3)\Rightarrow (1)$.

Assume that $R$ is a right coherent ring. Notice that the
subcategory of $\Mod R$ consisting of G-flat modules is closed under
extensions and direct summands by [15, Theorem 3.7], and also notice that a
G-flat module in $\Mod R$ with finite flat dimension is flat by
Lemma 2.3, so we get $(5)\Rightarrow (1)$ by using an argument
similar to that in the proof of $(5)\Rightarrow (1)$ in Theorem
3.9. $\hfill{\square}$

\vspace{0.2cm}

From the above argument, we see that $n$-SG-Flat$(R)$ is not closed
under direct summands in general. However, the following result,
which is a generalization of [18, Lemma 2.3], shows that some
special direct summand of an $n$-SG-flat module is again
$n$-SG-flat. For a module $M\in \Mod R$, we use $\widetilde{M}$ to
denote the maximal submodule of $M$ without flat direct summands.

\vspace{0.2cm}

{\bf Theorem 4.5} {\it For any $n\geq 1$, a module $M\in \Mod R$ is
$n$-SG-flat if and only if so is $\widetilde{M}$.}

\vspace{0.2cm}

{\it Proof.} The sufficiency follows from Proposition 4.3. In the
following, we will prove the necessity.

Assume that $M\in \Mod R$ is $n$-SG-flat and $I\in \Mod R^{op}$ is
any injective module. Then there exists an exact sequence:
$$0\rightarrow (M=)\widetilde{M}\oplus F\stackrel{h_{n}}{\longrightarrow} F_{n-1}
\stackrel{h_{n-1}}{\longrightarrow}
\cdots\stackrel{h_{1}}{\longrightarrow}F_{0}\stackrel{h_{0}}{\longrightarrow}\widetilde{M}\oplus
F(=M) \rightarrow 0$$ in $\Mod R$ with $F$ and $F_{i}$ flat for any
$0\leq i\leq n-1$, such that $I\otimes_{R}-$ leaves the sequence
exact. Put $\Im h_{i}=L_{i}$ for any $0\leq i \leq n$. We first
consider the following push-out diagram:
$$\xymatrix@=18pt
{&&0\ar[d]&0\ar[d]\\
0\ar[r]&F\ar[r]\ar@{=}[d]&M\ar[d]^{h_n}\ar[r]&\widetilde{M}\ar[r]\ar[d]&0\\
0\ar[r]&F\ar[r]&F_{n-1}\ar[d]\ar[r]&H_{n-1}\ar[d]\ar[r]&0\\
&&L_{n-1}\ar[d]\ar@{=}[r]&L_{n-1}\ar[d]\\
&&0&0 } $$ Then we have the following diagram with exact columns and
rows:
$$\xymatrix@C=20pt@R=15pt {
&0\ar[d]&0\ar[d]\\
&L_{n-1}^+\ar[d]\ar@{=}[r]&L_{n-1}^+\ar[d]\\
0\ar[r]&H_{n-1}^+\ar[r]\ar[d]&F_{n-1}^+\ar[d]\ar[r]&F^+\ar[r]\ar@{=}[d]&0\\
0\ar[r]&(\widetilde{M})^+\ar[d]\ar[r]&M^+\ar[d]\ar[r]&F^+\ar[r]&0\\
&0&0}\ \ \
$$
where both $F^+$ and $F_{n-1}^+$ are injective by [12, Theorem 2.1].
Because both $M^+$ and $L_{n-1}^+$ are G-injective by [15, Theorem 3.6],
both $(\widetilde{M})^+$ and $H_{n-1}^+$ are also G-injective by [15, Theorem 2.6].
Thus $\Ext_{R}^1(F^+, H_{n-1}^+)=0$ and the middle row
$$0\rightarrow H_{n-1}^+\rightarrow F_{n-1}^+\rightarrow
F^+\rightarrow 0$$ in the above diagram splits. So $H_{n-1}^+$ is
injective and hence $H_{n-1}$ is flat again by [12, Theorem 2.1]. Because
$L_{n-1}$ is $n$-SG-flat, the third column
$$0\rightarrow \widetilde{M}\rightarrow H_{n-1}\rightarrow
L_{n-1}\rightarrow 0$$ in the former diagram is still exact after applying
the functor $I\otimes_{R} -$.

Next, we consider the following pull-back diagram:
$$\xymatrix@C=20pt@R=15pt
{&0\ar[d]&0\ar[d]\\
&L_1\ar[d]\ar@{=}[r]&L_1\ar[d]\\
0\ar[r]&H_0\ar[r]\ar[d]&F_0\ar[d]\ar[r]&F\ar[r]\ar@{=}[d]&0\\
0\ar[r]&\widetilde{M}\ar[d]\ar[r]&M\ar[d]\ar[r]&F\ar[r]&0\\
&0&0}\ \ \
$$
From the middle row in the above diagram we know that $H_0$ is flat.
Then from the third row in the above diagram, we get an exact
sequence: $$0=\Tor _{i+1}^{R}(I, F)\rightarrow \Tor _{i}^{R}(I,
\widetilde{M})\rightarrow \Tor _{i}^{R}(I, M)=0$$ for any $i\geq 1$.
Thus $\Tor_{i}^{R}(I, \widetilde{M})=0$ for any $i\geq 1$ and
therefore $0\rightarrow I\otimes_{R}L_1\rightarrow
I\otimes_{R}H_0\rightarrow I\otimes_{R}\widetilde{M}\rightarrow 0$
is exact. So we obtain the following exact sequence:
$$0\rightarrow \widetilde{M}\rightarrow H_{n-1}\rightarrow
F_{n-2}\rightarrow\cdots\rightarrow F_{1}\rightarrow H_0\rightarrow
\widetilde{M}\rightarrow 0,$$which is still exact after applying the
functor $I\otimes_{R} -$, which implies that $\widetilde{M}$ is
$n$-SG-flat. $\hfill{\square}$

\vspace{0.2cm}

We call two modules $M, N\in \Mod R$ {\it flatly equivalent}
if there exist flat modules $F_1, F_2$ in $\Mod R$, such that
$M\oplus F_1 \cong N\oplus F_2$. By Theorems 4.5, we immediately
have the following

\vspace{0.2cm}

{\bf Corollary 4.6} {\it Assume that $M, N\in \Mod R$ are flatly
equivalent. Then, for any $n\geq 1$, $M$ is $n$-SG-flat if and only
if so is $N$.}

\vspace{0.2cm}

In the rest of this section, we will investigate the relation
between $n$-SG-flat modules and $n$-SG-projective (resp. injective)
modules. We first have the following result, which is a
generalization of [3, Proposition 3.9].

\vspace{0.2cm}

{\bf Proposition 4.7} {\it For any $n\geq 1$, a finitely generated
$n$-SG-projective $R$-module is finitely presented $n$-SG-flat.}

\vspace{0.2cm}

{\it Proof.} Assume that $M$ is a finitely generated
$n$-SG-projective $R$-module. By Theorem 3.14, there exists an exact
sequence:
$$0\rightarrow M\rightarrow P_{n-1}\rightarrow
P_{n-2}\rightarrow\cdots\rightarrow P_0\rightarrow M\rightarrow 0$$
in $\mod R$ with $P_i$ projective for any $0\leq i\leq n-1$, and
$\Ext ^{i}_{R}(M, R) = 0$ for any $i\geq 1$. Let $I\in \Mod R^{op}$
be injective. By [5, Chapter VI, Proposition 5.3 ``Remark"], we have an isomorphism:
$$\Tor_ {i}^{R}(I, M)\cong \Hom _{R}(\Ext ^{i}_{R}(M,
R), I)$$ for any $i\geq 1$. Thus $\Tor _{i}^{R}(I, M) = 0$ for any
$i\geq 1$, and therefore $M$ is finitely presented $n$-SG-flat.
$\hfill{\square}$

\vspace{0.2cm}

As an application of Proposition 4.7, we give another example of
2-SG-flat modules, but not 1-SG-flat.

\vspace{0.2cm}

{\bf Example 4.8} Consider a Noetherian local ring $R=
k[[X,Y]]/(XY)$, where $k$ is a field. Then the ideals $(X+(XY))$ and
$(Y+(XY))$ of $R$ are finitely generated 2-SG-projective $R$-modules
by [6, Example 4.15], where $(X+(XY))$ and $(Y+(XY))$ are the
residue classes in $R$ of $X$ and $Y$. By Proposition 4.7, both
$(X+(XY))$ and $(Y+(XY))$ are 2-SG-flat, but neither of
them are 1-SG-flat by [3, Example 3.11].

\vspace{0.2cm}

The following result generalizes [18, Theorems 2.4 and 2.12].

\vspace{0.2cm}

{\bf Proposition 4.9} {\it (1) If $M\in \Mod R$ is $n$-SG-flat, then
$M^+\in \Mod R^{op}$ is $n$-SG-injective.

(2) For an Artinian algebra $R$, if $M\in \Mod R$ is
$n$-SG-injective, then $M^+\in \Mod R^{op}$ is $n$-SG-flat.}

\vspace{0.2cm}

{\it Proof.} (1) Assume that $M\in \Mod R$ is $n$-SG-flat. Then
there exists an exact sequence:
$$0\rightarrow M\rightarrow F_{n-1}\rightarrow
F_{n-2}\rightarrow\cdots\rightarrow F_0\rightarrow M\rightarrow 0$$
in $\Mod R$ with $F_{i}$ flat for any $0 \leq i \leq n-1$, and $\Tor
^{R}_{i}(E,M) = 0$ for any injective right $R$-module $E$ and $i\geq
1$. So we get the following exact sequence: $$0\rightarrow
M^+\rightarrow F^{+}_{0}\rightarrow\cdots\rightarrow
F^{+}_{n-2}\rightarrow F^{+}_{n-1}\rightarrow M^+\rightarrow 0$$ in
$\Mod R^{op}$ with $F^{+}_{i}$ injective by [12, Theorem 2.1]. Because
$\Ext ^{i}_{R}(E, M^{+})\cong \Tor ^{R}_{i}(E, M)^{+}= 0$ for any
$i\geq 1$ by [5, Chapter VI, Proposition 5.1], $M^+$ is $n$-SG-injective.

(2) Assume that $M\in \Mod R$ is $n$-SG-injective. Then there exists
an exact sequence:
$$0\rightarrow M\rightarrow
I^{n-1}\rightarrow I^{n-2}\rightarrow\cdots\rightarrow
I^0\rightarrow M\rightarrow 0$$ in $\Mod R$ with $I^{i}$ injective
for any $0 \leq i\leq n-1$. So we get the following exact sequence:
$$0\rightarrow M^+\rightarrow (I^{0})^{+}\rightarrow\cdots\rightarrow
(I^{n-2})^{+}\rightarrow (I^{n-1})^{+}\rightarrow M^+\rightarrow 0$$
in $\Mod R^{op}$ with $(I^{i})^{+}$ flat for any $0\leq i \leq n-1$.

Let $E\in \Mod R$ be any injective module. Then $E=\oplus_{\gamma\in
\Gamma}E_{\gamma}$ with $E_{\gamma} \in \mod R$ injective. By
[5, Chapter VI, Proposition 5.3 ``Remark"], we have the following isomorphism:
$$\Tor ^{R}_{i}(M^{+}, E)\cong\Tor ^{R}_{i}(M^{+}, \oplus_{\gamma
\in \Gamma} E_{\gamma})\cong\oplus_{\gamma \in \Gamma} \Tor
^{R}_{i}(M^{+}, E_{\gamma})\cong\oplus_{\gamma \in\Gamma} (\Ext
^{i}_{R}(E_{\gamma}, M))^{+}=0$$ for any $ i\geq 1$, which implies
that $M^+$ is $n$-SG-flat. $\hfill{\square}$

\vspace{0.5cm}

{\bf Acknowledgements} The research was partially supported by the
Specialized Research Fund for the Doctoral Program of Higher
Education (Grant No. 20060284002), NSFC (Grant No. 10771095) and NSF
of Jiangsu Province of China (Grant No. BK2007517). The authors
thank Dr. Xiaojin Zhang for his useful suggestions.

\end{document}